\documentclass[12pt]{article}
\usepackage{amsmath,amsfonts,epsfig}
\title{\bf Quotient Isomorphism Invariants of a Finitely Generated Coxeter Group}
\author{Michael Mihalik, John Ratcliffe, and Steven Tschantz \\ 
Mathematics Department, Vanderbilt University, \\
Nashville TN 37240, USA}
\newtheorem{theorem}{Theorem}[section]
\newtheorem{proposition}[theorem]{Proposition}

\newtheorem{corollary}[theorem]{Corollary}

\newenvironment{proof}{{\bf Proof:\ }}{\hfill$\square$\vspace{.2in}}

\def\ov{\overline}

\date{}
\begin{document}
\maketitle

\section{Introduction} 

The isomorphism problem for finitely generated Coxeter groups is the problem of deciding 
if two finite Coxeter matrices define 
isomorphic Coxeter groups. Coxeter \cite{Coxeter} solved this problem for finite irreducible 
Coxeter groups.  
Recently there has been considerable interest and activity on the isomorphism problem 
for arbitrary finitely generated Coxeter groups. 

In this paper we describe a family of isomorphism invariants 
of a finitely generated Coxeter group $W$.  
Each of these invariants is the isomorphism type 
of a quotient group $W/N$ of $W$ by a characteristic subgroup $N$. 
The virtue of these invariants is that $W/N$ is also a Coxeter group. 
For some of these invariants, the isomorphism problem of $W/N$ is solved 
and so we obtain isomorphism invariants that can be effectively used 
to distinguish isomorphism types of finitely generated Coxeter groups. 

We emphasize that even if the isomorphism problem for 
finitely generated Coxeter groups is eventually solved, 
several of the algorithms described in our paper will still be useful because 
they are computational fast and would most likely be incorporated into an efficient 
computer program that determines if two finite 
rank Coxeter systems have isomorphic groups.  

In \S2, we establish notation. 
In \S3, we describe two elementary quotienting operations on a Coxeter system 
that yields another Coxeter system. 
In \S4, we describe the binary isomorphism invariant of a finitely generated Coxeter group. 
In \S5, we review some matching theorems. 
In \S6, we describe the even isomorphism invariant of a finitely generated Coxeter group. 
In \S7, we define basic characteristic subgroups of a finitely generated Coxeter group.   
In \S8, we describe the spherical rank two isomorphism invariant of a finitely generated 
Coxeter group. 
In \S9, we make some concluding remarks.

\section{Preliminaries} 
 
A {\it Coxeter matrix} is a symmetric matrix $M = (m(s,t))_{s,t\in S}$ 
with $m(s,t)$ either a positive integer or infinity and $m(s,t) =1$ 
if and only if $s=t$. A {\it Coxeter system} with Coxeter matrix $M = (m(s,t))_{s,t\in S}$ 
is a pair $(W,S)$ consisting of a group $W$ and a set of generators $S$ for $W$ 
such that $W$ has the presentation
$$W =\langle S \ |\ (st)^{m(s,t)}:\, s,t \in S\ \hbox{and}\ m(s,t)<\infty\rangle.$$
We call the above presentation of $W$, the {\it Coxeter presentation} of $(W,S)$. 
If $(W,S)$ is a Coxeter system with Coxeter matrix $M = (m(s,t))_{s,t\in S}$, 
then the order of $st$ is $m(s,t)$ for each $s,t$ in $S$ by 
Prop. 4, p. 92 of Bourbaki \cite{Bourbaki}, and so a Coxeter system $(W,S)$ 
determines its Coxeter matrix; moreover, any Coxeter matrix 
$M = (m(s,t))_{s,t\in S}$ determines a Coxeter system $(W,S)$ where $W$ 
is defined by the corresponding Coxeter presentation. 
If $(W,S)$ is a Coxeter system, 
then $W$ is called a {\it Coxeter group} and $S$ is called a set of {\it Coxeter generators} 
for $W$, and the cardinality of $S$ is called the {\it rank} of $(W,S)$.
A Coxeter system $(W,S)$ has finite rank if and only if $W$ 
is finitely generated by Theorem 2 (iii), p. 20 of Bourbaki  \cite{Bourbaki}. 

Let $(W,S)$ be a Coxeter system. A {\it visual subgroup} of $(W,S)$ 
is a subgroup of $W$ of the form $\langle A\rangle$ for some $A \subset S$. 
If $\langle A\rangle$ is a visual subgroup of $(W,S)$, 
then $(\langle A\rangle, A)$ is also a Coxeter system 
by Theorem 2 (i), p. 20 of Bourbaki  \cite{Bourbaki}.

When studying a Coxeter system $(W,S)$ with Coxeter matrix $M$ 
it is helpful to have a visual representation of $(W,S)$. 
There are two graphical ways of representing $(W,S)$ 
and we shall use both depending on our needs. 

The {\it Coxeter diagram} ({\it {\rm C}-diagram}) of $(W,S)$ is the labeled undirected graph 
$\Delta = \Delta(W,S)$ with vertices $S$ and edges 
$$\{(s,t) : s, t \in S\ \hbox{and}\ m(s,t) > 2\}$$
such that an edge $(s,t)$ is labeled by $m(s,t)$. 
Coxeter diagrams are useful for visually representing finite Coxeter groups. 
If $A\subset S$, then $\Delta(\langle A\rangle,A)$ is the 
subdiagram of $\Delta(W,S)$ induced by $A$.

A Coxeter system $(W,S)$ is said to be {\it irreducible} 
if its C-diagram $\Delta$ is connected. 
A visual subgroup $\langle A\rangle$ of $(W,S)$ is said to be {\it irreducible} 
if $(\langle A\rangle, A)$ is irreducible. 
A subset $A$ of $S$ is said to be {\it irreducible} if $\langle A\rangle$ is irreducible. 

A subset $A$ of $S$ is said to be a {\it component} of $S$ if $A$ is a maximal irreducible 
subset of $S$ or equivalently if $\Delta(\langle A\rangle, A)$ is a connected component 
of $\Delta(W,S)$. 
The connected components of the $\Delta(W,S)$
represent the factors of a direct product decomposition of $W$.

The {\it presentation diagram} ({\it {\rm P}-diagram}) of $(W,S)$ is the labeled undirected graph 
$\Gamma = \Gamma(W,S)$ with vertices $S$ and edges 
$$\{(s,t) : s, t \in S\ \hbox{and}\ m(s,t) < \infty\}$$
such that an edge $(s,t)$ is labeled by $m(s,t)$. 
Presentation diagrams are useful for visually representing infinite Coxeter groups. 
If $A\subset S$, then $\Gamma(\langle A\rangle,A)$ is the 
subdiagram of $\Gamma(W,S)$ induced by $A$. 
The connected components of $\Gamma(W,S)$
represent the factors of a free product decomposition of $W$. 

For example, consider the Coxeter group $W$ generated by the four reflections 
in the sides of a rectangle in $E^2$.  The C-diagram of $(W,S)$ is 
the disjoint union of two edges labeled by $\infty$ while  
the P-diagram of $W$ is a square with edge labels 2.   

Let $(W,S)$ and $(W',S')$ be Coxeter systems 
with P-diagrams $\Gamma$ and $\Gamma'$, respectively. 
An {\it isomorphism} $\phi: (W,S) \to (W',S')$ of Coxeter systems 
is an isomorphism $\phi: W\to W'$ such that $\phi(S) = S'$. 
An {\it isomorphism} $\psi:\Gamma\to \Gamma'$ of P-diagrams is a bijection 
from $S$ to $S'$ that preserves edges and their labels. 
Note that $(W,S) \cong (W',S')$ if and only if $\Gamma \cong \Gamma'$. 

We shall use Coxeter's notation on p. 297 of \cite{Coxeterb}
for the irreducible spherical Coxeter simplex reflection groups except 
that we denote the dihedral group ${\bf D}_2^k$ by ${\bf D}_2(k)$. 
Subscripts denote the rank of a Coxeter system in Coxeter's notation. 
Coxeter's notation partly agrees with but differs from Bourbaki's notation on p.193 of 
\cite{Bourbaki}.

Coxeter \cite{Coxeter} proved that every finite irreducible Coxeter system 
is isomorphic to exactly one 
of the Coxeter systems ${\bf A}_n$, $n\geq 1$, ${\bf B}_n$, $n\geq 4$, ${\bf C}_n$, 
$n\geq 2$, ${\bf D}_2(k)$, $k\geq 5$, 
${\bf E}_6$, ${\bf E}_7$, ${\bf E}_8$, ${\bf F}_4$, ${\bf G}_3$, ${\bf G}_4$. 
For notational convenience, we define ${\bf B}_3 = {\bf A}_3$, ${\bf D}_2(3) = {\bf A}_2$, 
and ${\bf D}_2(4) = {\bf C}_2$

The {\it type} of a finite irreducible Coxeter system $(W,S)$ is the 
isomorphism type of $(W,S)$ represented by one of the systems 
${\bf A}_n$, ${\bf B}_n$, ${\bf C}_n$, ${\bf D}_2(k)$, 
${\bf E}_6$, ${\bf E}_7$, ${\bf E}_8$, ${\bf F}_4$, ${\bf G}_3$, ${\bf G}_4$. 
The {\it type} of an irreducible subset $A$ of $S$  
is the type of $(\langle A\rangle,A)$.  

The C-diagram of ${\bf A}_n$ is a linear diagram 
with $n$ vertices and all edge labels 3. 
The C-diagram of ${\bf B}_n$ is a Y-shaped diagram 
with $n$ vertices and all edge labels 3 and two short 
arms of consisting of single edges. 
The C-diagram of ${\bf C}_n$ is a linear diagram with $n$ vertices 
and all edge labels 3 except for the last edge labelled 4. 
The C-diagram of ${\bf D}_2(k)$ is a single edge with label $k$. 
The C-diagrams of ${\bf E}_6$, ${\bf E}_7, {\bf E}_8$ 
are star shaped with three arms and all edge labels 3. 
One arm has length one and another has length two. 
The C-diagram of ${\bf F}_4$ is a linear diagram 
with edge labels $3,4,3$ in that order.  
The C-diagram of ${\bf G}_3$ is a linear diagram 
with edge labels $3, 5$. 
The C-diagram of ${\bf G}_4$ is a linear diagram 
with edge labels $3,3,5$ in that order. 

\section{Elementary Quotient Operations}  

In this section we describe two types of elementary edge quotient 
operations on a Coxeter system $(W,S)$ of finite rank. 
The first we call edge label reduction and  
the second we call edge elimination. 
 
Suppose $s$ and $t$ are distinct elements of $S$ with $m(s,t) < \infty$. 
Let $d$ be a positive divisor of $m = m(s,t)$, with $d < m$, 
and let $N$ be the normal closure of the element $(st)^d$ of $W$. 
Then a presentation for $W/N$ is obtained 
from the Coxeter presentation for $(W,S)$ by adding 
the relator $(st)^d$.  As $m = (m/d)d$, 
the relator $(st)^m$ is derivable from the relator $(st)^d$ 
and so the relator $(st)^m$ can be removed from the presentation for $W/N$. 

Assume $d > 1$. 
Then the presentation for $W/N$ is a Coxeter presentation whose 
P-diagram is obtained from the P-diagram for $(W,S)$ 
by replacing the label $m$ on the edge $(s,t)$ 
with the label $d$. 
We call the operation of passing from 
the Coxeter system $(W,S)$ to the quotient Coxeter system 
$(W/N, \{sN: s\in S\})$ the $(s,t)$ {\it edge label reduction} 
from $m$ to $d$. 
For example, if we reduce the 4 edge of ${\bf F}_4$ to 2,  
we obtain the Coxeter system ${\bf A}_2 \times {\bf A}_2$. 

Now assume $d=1$. 
We delete from the presentation 
for $W/N$ the generator $t$ 
and the relator $st$ and replace all occurrences of $t$ 
in the remaining relators by $s$. 
Suppose $r$ is in $S$ and $k = m(r,s) < \infty$ 
and $\ell = m(r,t) < \infty$.  
Then we have the relators $(rs)^k$ and $(rs)^\ell$ in the presentation for $W/N$. 
Let $d$ be the greatest common divisor of $k$ and $\ell$. 
Then there are integers $a$ and $b$ such that $d = ak+b\ell$. 
This implies that $(rs)^d$ is derivable from $(rs)^k$ and $(rs)^\ell$ 
and so we may add the relator $(rs)^d$ to the presentation for $W/N$. 
Then $(rs)^k$ and $(rs)^\ell$ are derivable from $(rs)^d$ 
and so we can eliminate the relators $(rs)^k$ and $(rs)^\ell$ 
from the presentation for $W/N$. 
We do this for each $r$ in $S$ 
such that $m(r,s) < \infty$ and $m(r,t) < \infty$. 
On the $P$-diagram level, we have eliminated the edge $(s,t)$ 
and identified the vertices $s$ and $t$ and we have coalesced 
each edge $(r,s)$ with label $k< \infty$ with the edge $(r,t)$ with label $\ell < \infty$ 
to form an edge with label $d$ the greatest common divisor of $k$ and $\ell$. 
If each common divisor $d$ is greater than one, we obtain a Coxeter 
presentation for $W/N$.  If some common divisor $d$ is one, 
we delete the corresponding generator $r$ and repeating the above reduction procedure 
on the presentation of $W/N$. 
As the set $S$ of generators is finite, we will eventually stop deleting generators 
and obtain a Coxeter presentation for $W/N$ 
with generators the subset $S'$ of $\{sN: s\in S\}$ 
corresponding to the undeleted elements of $S$. 
We call the operation of passing from 
the Coxeter system $(W,S)$ to the quotient Coxeter system 
$(W/N, S')$ the $(s,t)$ {\it edge elimination}. 
For example, if we eliminate the 3 edge from ${\bf C}_3$, 
we obtain the Coxeter system ${\bf A}_1\times {\bf A}_1$. 

\section{The Binary Isomorphism Invariant}

Let $(W,S)$ be a Coxeter system of finite rank. 
For each pair of elements $s,t$ of $S$ with $m(s,t) < \infty$, 
let $b(s,t)$ be the 2-part of $m(s,t)$, that is, $b(s,t)$ is the largest 
power of 2 that divides $m(s,t)$. 
If $m(s,t) = \infty$, we define $b(s,t) = \infty$. 
Let $N_b$ be the normal closure in $W$ 
of all the elements of the form $(st)^{b(s,t)}$ with $b(s,t) < \infty$, 
and let $W_b = W/N_b$. 
Let $\eta:W \to W_b$ be the quotient homomorphism, 
and let $S_b=\eta(S)$.  

\begin{proposition} 
The pair $(W_b,S_b)$ is a Coxeter system such that if $s$ and $t$ 
are in $S$, then $\eta(s) = \eta(t)$ 
if and only if $s$ and $t$ are conjugate in $W$.  
If $s$ and $t$ are nonconjugate elements of $S$, 
then the order of $\eta(s)\eta(t)$ is the minimum of the set of all $b(u,v)$ 
such that $u$ and $v$ are in $S$ 
and $u$ is conjugate to $s$ and $v$ is conjugate to $t$. 
In particular, the order of $\eta(s)\eta(t)$ is a power of 2 or $\infty$. 
\end{proposition}

\noindent\begin{proof}
The system $(W_b,S_b)$ can be obtained from $(W,S)$ by a sequence of 
elementary quotient operations. First we can do a series of edge label 
reductions of all the even labelled edges of the P-diagram of $(W,S)$ 
to their 2-parts. 
Then we do a series of edge eliminations of all the odd labelled edges.  
Each element of the form $(st)^{b(s,t)}$ with $b(s,t) < \infty$ 
is in the commutator subgroup of $W$. 
Therefore abelianizing $W$ factors through the quotient $W/N_b$, 
and so $\eta(s) = \eta(t)$ if and only if $s$ and $t$ are the same odd component 
of the P-diagram of $(W,S)$.  Hence $\eta(s) = \eta(t)$ if and only if 
$s$ and $t$ are conjugate in $W$ by Prop. 3, p. 12 of Bourbaki \cite{Bourbaki}. 

Suppose $s$ and $t$ are nonconjugate elements of $S$ 
and $u$ and $v$ are in $S$, with $m(u,v) < \infty$,  
and $u$ is conjugate to $s$ and $v$ is conjugate to $t$. 
Then $u$ and $v$ are not conjugate, 
and so $m(u,v)$ is even, and therefore $b(u,v)$ is 
a power of 2 greater than 1. 
In the coalescence of two such edges, the greatest common divisor is the 
minimum of the two edge labels. Therefore 
the order of $\eta(s)\eta(t)$ is the minimum 
of the set of all $b(u,v)$ such that $u$ and $v$ are in $S$ 
and $u$ is conjugate to $s$ and $v$ is conjugate to $t$. 
\end{proof}

\begin{theorem} 
Let $(W,S)$ be a Coxeter system of finite rank. 
For each pair of elements $s,t$ of $S$ with $m(s,t) < \infty$, 
let $b(s,t)$ be the largest power of 2 that divides $m(s,t)$. 
Let $N_b$ be the normal closure in $W$ 
of all the elements of the form $(st)^{b(s,t)}$ with $m(s,t) < \infty$, 
Then $N_b$ is the normal closure in $W$ 
of the set of all elements of $W$ of odd order. 
Therefore $N_b$ is a characteristic subgroup of $W$ 
that does not depend on the choice of Coxeter generators $S$.  
\end{theorem}
\noindent\begin{proof} 
Every element of the form $(st)^{b(s,t)}$ with $m(s,t) < \infty$ 
has odd order, and so $N_b$ is contained in 
the normal closure of all the elements of odd order. 
Let $w$ be an element of odd order, 
then $\eta(w)$ has odd order in $W_b = W/N_b$. 
By the previous proposition, $(W_b,S_b)$ is a Coxeter system with all edge labels a power of 2.  
Therefore $\eta(w)$ is conjugate 
to an element of odd order of a finite 
visual subgroup of $(W_b,S_b)$ by \cite{Bourbaki}, Ch. V, \S 4, Ex. 2. 
The finite visual subgroups of $(W_b,S_b)$ are direct products of 
groups of type ${\bf A}_1$ and ${\bf C}_2$, and so are 2-groups. 
Therefore $W_b$ has no nontrivial elements of odd order. 
Hence $\eta(w) = 1$, and so $w$ is in $N_b$. 
Thus $N_b$ is the normal closure of all the elements of $W$ of odd order. 
\end{proof} 

P. Bahls proved in his Ph.D. thesis \cite{Bahls} 
that any finitely generated Coxeter group 
has at most one P-diagram, up to isomorphism, 
with all edge labels even;  
see Theorem 5.2 in Bahls and Mihalik \cite{B-M}. 
Therefore the isomorphism type of the P-diagram of $(W_b,S_b)$ 
is an isomorphism invariant of $W$ by Theorem 4.2.  
We call the isomorphism type of the P-diagram of $(W_b,S_b)$ 
the {\it binary isomorphism invariant} of $W$. 

In Figure 1, we illustrate two P-diagrams and 
their binary isomorphism invariant P-diagrams below them. 
The even diagrams are not isomorphic, and so 
the top two P-diagrams represent nonisomorphic Coxeter groups.

$$\mbox{
\setlength{\unitlength}{.8cm}
\begin{picture}(13,7)(0,0)
\thicklines
\put(.5,1){\circle*{.15}}
\put(.5,1){\line(1,0){1.5}}
\put(2,1){\circle*{.15}}
\put(2,1){\line(2,-1){1.5}}
\put(2,1){\line(2,1){1.5}}
\put(3.5,.25){\line(0,1){1.5}}
\put(3.5,.25){\circle*{.15}}
\put(3.5,1.75){\circle*{.15}}
\put(1.25,1.35){2}
\put(2.5,0){4}
\put(2.5,1.75){4}
\put(3.8,.9){2}
\put(0,4){\circle*{.15}}
\put(0,4){\line(1,0){2}}
\put(2,4){\circle*{.15}}
\put(0,4){\line(1,-1){1}}
\put(1,3){\circle*{.15}}
\put(2,4){\line(-1,-1){1}}
\put(2,4){\line(0,1){1.5}}
\put(2,5.5){\circle*{.15}}
\put(2,4){\line(2,-1){1.5}}
\put(3.5,3.25){\circle*{.15}}
\put(2,5.5){\line(2,1){1.5}}
\put(3.5,6.25){\circle*{.15}}
\put(3.5,6.25){\line(2,-1){1.5}}
\put(5,5.5){\circle*{.15}}
\put(3.5,3.25){\line(2,1){1.5}}
\put(5,4){\circle*{.15}}
\put(5,4){\line(0,1){1.5}}
\put(.1,3){2}
\put(1.6,3){3}
\put(2.5,3.1){3}
\put(4.2,3.1){12}
\put(.8,4.25){4}
\put(1.55,4.6){4}
\put(5.3,4.6){3}
\put(2.5,6.2){3}
\put(4.2,6.2){6}
\put(8.5,1){\circle*{.15}}
\put(8.5,1){\line(1,0){1.5}}
\put(10,1){\circle*{.15}}
\put(10,1){\line(2,-1){1.5}}
\put(10,1){\line(2,1){1.5}}
\put(11.5,.25){\line(0,1){1.5}}
\put(11.5,.25){\circle*{.15}}
\put(11.5,1.75){\circle*{.15}}
\put(9.25,1.35){2}
\put(10.5,0){2}
\put(10.5,1.75){4}
\put(11.8,.9){4}
\put(8,4){\circle*{.15}}
\put(8,4){\line(1,0){2}}
\put(10,4){\circle*{.15}}
\put(8,4){\line(1,-1){1}}
\put(9,3){\circle*{.15}}
\put(10,4){\line(-1,-1){1}}
\put(10,4){\line(0,1){1.5}}
\put(10,5.5){\circle*{.15}}
\put(10,4){\line(2,-1){1.5}}
\put(11.5,3.25){\circle*{.15}}
\put(10,5.5){\line(2,1){1.5}}
\put(11.5,6.25){\circle*{.15}}
\put(11.5,6.25){\line(2,-1){1.5}}
\put(13,5.5){\circle*{.15}}
\put(11.5,3.25){\line(2,1){1.5}}
\put(13,4){\circle*{.15}}
\put(13,4){\line(0,1){1.5}}
\put(8.1,3){2}
\put(9.6,3){3}
\put(10.5,3.1){3}
\put(12.2,3.1){6}
\put(8.8,4.25){4}
\put(9.55,4.6){4}
\put(13.3,4.6){3}
\put(10.5,6.2){3}
\put(12.2,6.2){12}
\end{picture}}$$

\medskip

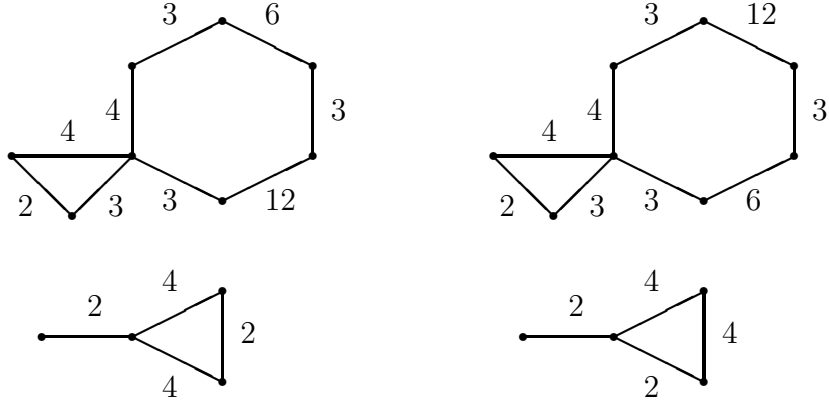
\centerline{\bf Figure 1}

\medskip

\section{Matching Theorems}

Let $(W,S)$ be a Coxeter system. 
A {\it basic subgroup} of $(W,S)$ is a noncyclic, maximal, finite, irreducible, 
visual subgroup of $(W,S)$. 
A {\it base} of $(W,S)$ is a subset $B$ of $S$ such that 
$\langle B\rangle$ is a basic subgroup of $(W,S)$. 
The theorems in this section are proved in our paper \cite{R-M-T}.

\begin{theorem} {\rm (Basic Matching Theorem)} 
Let $W$ be a finitely generated Coxeter group with 
two sets of Coxeter generators $S$ and $S'$. 
Let $B$ be a base of $(W,S)$. 
Then there is a unique irreducible subset $B'$ of $S'$ such that 
$[\langle B\rangle,\langle B\rangle]$ is conjugate to 
$[\langle B'\rangle,\langle B'\rangle]$ in $W$. Moreover, 

\begin{enumerate}
\item The set $B'$ is a base of $(W,S')$, 

\item If $|\langle B\rangle|=|\langle B'\rangle|$, then $B$ and $B'$ have the same type 
and there is an isomorphism $\phi:\langle B\rangle \to \langle B'\rangle$ 
that restricts to conjugation on $[\langle B\rangle,\langle B\rangle]$ 
by an element of $W$. 

\item If $|\langle B\rangle|<|\langle B'\rangle|$, then either
$B$ has type ${\bf B}_{2q+1}$ and 
$B'$ has type ${\bf C}_{2q+1}$ for some $q\geq 1$ or 
$B$ has type ${\bf D}_2(2q+1)$ and 
$B'$ has type ${\bf D}_2(4q+2)$ for some $q\geq 1$. 
Moreover, there is a monomorphism $\phi:\langle B\rangle \to \langle B'\rangle$ 
that restricts to conjugation on $[\langle B\rangle,\langle B\rangle]$ 
by an element of $W$. 
\end{enumerate}
\end{theorem}

Let $W$ be a finitely generated Coxeter group with two sets of Coxeter generators $S$ and $S'$. 
A basic subgroup $\langle B\rangle$ of $(W,S)$ is said to {\it match} 
a basic subgroup $\langle B'\rangle$ of $(W,S')$ 
if $[\langle B\rangle,\langle B\rangle]$ is conjugate to 
$[\langle B'\rangle,\langle B'\rangle]$ in $W$. 
A base $B$ of $(W,S)$ is said to {\it match} a base $B'$ of $(W,S')$ 
if $\langle B\rangle$ matches $\langle B'\rangle$.

\begin{theorem} 
Let $(W,S)$ be a Coxeter system of finite rank.  
Let $B$ be a base of $(W,S)$ of type ${\bf C}_{2q+1}$ for some $q\geq 1$, and   
let $a,b,c$ be the elements of $B$ such that $m(a,b)=4$ and $m(b,c)=3$. 
Then $W$ has a set of Coxeter generators $S'$ such that $B$ matches a base $B'$ 
of $(W,S')$ of type ${\bf B}_{2q+1}$ if and only if 
$m(s,t)=2$ for all $(s,t)\in (S-B)\times B$ such that $m(s,a)<\infty$. 
\end{theorem} 

\begin{theorem} 
Let $B$ be a base of $(W,S)$ of type ${\bf C}_{2q+1}$ for some $q\geq 1$, and   
let $a,b,c$ be the elements of $B$ such that $m(a,b)=4$ and $m(b,c)=3$. 
Suppose that $m(s,t)=2$ for all $(s,t)\in (S-B)\times B$ such that $m(s,a)<\infty$. 
Let $d =aba$, and let $z$ be the longest element of $\langle B\rangle$. 
Let $S'=(S-\{a\})\cup\{d,z\}$ and $B'=(B-\{a\})\cup\{d\}$. 
Then $S'$ is a set of Coxeter generators for $W$ 
such that 
\begin{enumerate}
\item The set $B'$ is a base of $(W,S')$ of type ${\bf B}_{2q+1}$ that matches $B$,  
\item $m(z,t)=2$ for all $t\in B'$,
\item If $(s,t)\in (S-B)\times\{d,z\}$, then 
$m(s,t)<\infty$ if and only if $m(s,a)<\infty$, moreover  
if $m(s,t)<\infty$, then $m(s,t)=2$.
\end{enumerate}
\end{theorem}

\begin{theorem} 
Let $(W,S)$ be a Coxeter system of finite rank, and   
let $B=\{a,b\}$ be a base of $(W,S)$ of type ${\bf D}_2(4q+2)$ for some $q\geq 1$. 
Then $W$ has a set of Coxeter generators $S'$ such that $B$ 
matches a base $B'$ of $(W,S')$ of type ${\bf D}_2(2q+1)$ if and only if 
either $v=a$ or $v=b$ has the property that 
if $s\in S-B$ and $m(s,v)<\infty$, then $m(s,a) = m(s,b) = 2$.  
\end{theorem}

\begin{theorem} 
Let $B=\{a,b\}$ be a base of $(W,S)$ of type ${\bf D}_2(4q+2)$ for some $q\geq 1$. 
Suppose that if $s\in S-B$ and $m(s,a)<\infty$, then $m(s,a) = m(s,b) = 2$. 
Let $c=aba$ and let $z$ be the longest element of $\langle B\rangle$. 
Let $S'=(S-\{a\})\cup\{c,z\}$ and $B'=\{b,c\}$. 
Then $S'$ is a set of Coxeter generators of $W$ such that 
\begin{enumerate}
\item The set $B'$ is a base of $(W,S')$ of type ${\bf D}_2(2q+1)$ that matches $B$,  
\item $m(z,b)=m(z,c) = 2$, 
\item if $(s,t)\in (S-B)\times \{c,z\}$, then $m(s,t)<\infty$ if and only if $m(s,a)<\infty$,  
moreover if $m(s,t)<\infty$, then $m(s,t) = 2$. 
\end{enumerate}
\end{theorem}

\section{The Even Isomorphism Invariant}  

Let $(W,S)$ be a Coxeter system of finite rank.  
A pair of elements $(a,b)$ of $S$ is said to be {\it special} 
if $(a,b)$ satisfy the conditions of Theorem 5.4, that is,  
if $m(a,b) \equiv 2\ {\rm mod}\ 4$ 
and either $v = a$ or $v=b$ has the property 
that if $s \in S-\{a,b\}$ and $m(s,v)<\infty$, then $m(s,a) = m(s,b) = 2$. 

Let $(a,b)$ be a pair of elements of $S$. 
If $(a,b)$ is special, define $\ov m(a,b)= 2$ 
and if $(a,b)$ is nonspecial, define $\ov m(a,b) = m(a,b)$.

Let $N_{e}$ be the normal closure in $W$ 
of all the elements of the form $ab$ with 
$a$ and $b$ elements of $S$ such that $m(a,b)$ is odd 
together with all the elements of the form $(ab)^2$ 
such that $(a,b)$ is a special pair of elements of $S$. 
Let $W_{e} = W/N_{e}$,  
let $\eta:W \to W_{e}$ be the quotient homomorphism, 
and let $S_{e}=\eta(S)$.  

\begin{proposition} 
The pair $(W_{e},S_{e})$ is a Coxeter system such that if $s$ and $t$ are in $S$, 
then $\eta(s) = \eta(t)$ if and only if $s$ and $t$ are conjugate in $W$. 
If $s$ and $t$ are nonconjugate elements of $S$, 
then the order of $\eta(s)\eta(t)$ is the greatest common divisor of the set of all $\ov m(u,v)$ 
such that $u$ and $v$ are in $S$ 
and $u$ is conjugate to $s$ and $v$ is conjugate to $t$. 
In particular, the order of $\eta(s)\eta(t)$ is either even or $\infty$. 
\end{proposition}

\noindent\begin{proof}
The system $(W_{e},S_{e})$ can be obtained from $(W,S)$ by 
a sequence of elementary quotient operations. 
First we reduce to 2 all the edge labels of special edges $(a,b)$ 
to obtain a Coxeter system with Coxeter matrix $\ov M = (\ov m(s,t))_{s,t\in S}$. 
Then we eliminate all the odd labelled edges.  
Each element of the form either $st$, with $m(s,t)$ odd,  
or $(st)^2$, with $(s,t)$ special, is in the commutator subgroup of $W$. 
Therefore abelianizing $W$ factors through the quotient $W/N_{e}$, 
and so $\eta(s) = \eta(t)$ if and only if $s$ and $t$ are in the same odd component 
of the P-diagram of $(W,S)$.  Hence $\eta(s) = \eta(t)$ if and only if 
$s$ and $t$ are conjugate in $W$ by Prop. 3, p. 12 of Bourbaki \cite{Bourbaki}. 

Suppose $s$ and $t$ are nonconjugate elements of $S$ 
and $u$ and $v$ are in $S$, with $\ov m(u,v) < \infty$,  
and $u$ is conjugate to $s$ and $v$ is conjugate to $t$. 
Then $u$ and $v$ are not conjugate, 
and so $\ov m(u,v)$ is even. 
In the coalescence of two such edges, the greatest common divisor is even. 
Therefore the order of $\eta(s)\eta(t)$ is the greatest common divisor  
of the set of all $\ov m(u,v)$ such that $u$ and $v$ are in $S$ 
and $u$ is conjugate to $s$ and $v$ is conjugate to $t$. 
\end{proof}

\newpage

Let $(W,S)$ be a Coxeter system of finite rank. 
A base $B$ of $(W,S)$ is said to be of {\it odd type} 
if there are elements $a$ and $b$ in $B$, with $m(a,b)$ odd. 
A base $B$ of $(W,S)$ is said to be {\it special} 
if $B=\{a,b\}$ with $(a,b)$ special.

\begin{theorem} 
Let $(W,S)$ be a Coxeter system of finite rank. 
Then $N_{e}$ is a characteristic subgroup of $W$ 
that does not depend on the choice of Coxeter generators $S$.  
\end{theorem}
\noindent\begin{proof}
Observe that $N_{e}$ is the normal closure in $W$ 
of the commutator subgroups of all the basic subgroups $\langle B\rangle$ of $(W,S)$ 
such that the base $B$ is either of odd type or special. 
Let $S'$ be another set of Coxeter generators of $W$. 
By the Basic Matching Theorem and Theorem 5.4, 
the group $N_{e}$ is also the normal closure in $W$ 
of the commutator subgroups of all the basic subgroups $\langle B'\rangle$ of $(W,S')$ 
such that the base $B'$ is either of odd type or special.  
Therefore $N_{e}$ is the normal closure in $W$ of all the elements 
of the form $s't'$ with $s'$ and $t'$ in $S'$ and $m(s',t')$ odd 
together with all the elements of the form $(s't')^2$ with $(s',t')$ 
a special pair of elements of $S'$.  
Thus $N_{e}$ is a characteristic subgroup of $W$ 
that does not depend on the choice of Coxeter generators $S$. 
\end{proof} 

Let $(W,S)$ be a Coxeter system of finite rank. 
P. Bahls proved in his Ph.D. thesis \cite{Bahls} 
that any finitely generated Coxeter group 
has at most one P-diagram, up to isomorphism, 
with all edge labels even. 
Therefore the isomorphism type of the P-diagram of $(W_{e},S_{e})$ 
is an isomorphism invariant of $W$ by Theorem 6.2.  
We call the isomorphism type of the P-diagram of $(W_{e},S_{e})$ 
the {\it even isomorphism invariant} of $W$. 
For example, the even isomorphism invariant of the 
system ${\bf D}_2(6)$ is 
the isomorphism type of the P-diagram 
of the system ${\bf A}_1\times {\bf A}_1$.

\section{Basic Characteristic Subgroups}  

Let $(W,S)$ be a Coxeter system of finite rank. 
A basic subgroup $\langle B\rangle$ of $(W,S)$ is said to be {\it reduced} 
if for every set of Coxeter generators $S'$ of $W$, 
the base $B$ matches a base $B'$ of $(W,S')$ 
with $|\langle B\rangle| \leq |\langle B'\rangle|$. 
Note that a basic subgroup $\langle B\rangle$ of $(W,S)$ is {\it unreduced} 
precisely when the base $B$ satisfies the conditions of either 
Theorem 5.2 or 5.4. 
In particular, every special base is unreduced. 

Let ${\cal F}$ be a family of finite irreducible 
Coxeter system isomorphism types.  
Let $N({\cal F})$ be the normal closure in $W$ 
of the commutator subgroups of all 
the reduced basic subgroups of $(W,S)$ of isomorphism type 
contained in ${\cal F}$ together with the commutator subgroups of 
all the unreduced basic subgroups of $(W,S)$ that match 
a reduced basic subgroup of another system $(W,S')$   
of isomorphism type contained in ${\cal F}$. 
Let $W({\cal F}) = W/N({\cal F})$,  
let $\eta: W \to W({\cal F})$ be the quotient homomorphism, 
and let $S({\cal F}) = \eta(S)$. 

\begin{theorem}  
The pair $(W({\cal F}), S({\cal F}))$ is a Coxeter system 
that can be obtained from $(W,S)$ be a finite series 
of elementary edge quotient operations. 
The group $N({\cal F})$ is a characteristic subgroup of $W$ 
that does not depend on the choice of Coxeter generators $S$.
\end{theorem}
\noindent\begin{proof}
Quotienting out the commutator subgroup of a basic subgroup $\langle B\rangle$ of $(W,S)$ 
can be realized by reducing all the even labelled edges of the 
C-diagram of $(\langle B\rangle, B)$ to 2 and eliminating all the 
odd labelled edges of the C-diagram. 
Therefore $(W({\cal F}), S({\cal F}))$ is a Coxeter system 
that can be obtained from $(W,S)$ be a finite series 
of elementary edge quotient operations. 

Let $S'$ be another set of Coxeter generators of $W$. 
By the Basic Matching Theorem,  
$N({\cal F})$ is also the normal closure in $W$ 
of the commutator subgroups of all the reduced basic subgroups of $(W,S')$ 
of isomorphism type contained in ${\cal F}$ together with the commutator subgroups 
of all the unreduced basic subgroups of $(W,S)$ that match 
a reduced basic subgroup of isomorphism type contained in ${\cal F}$. 
Thus $N({\cal F})$ is a characteristic subgroup of $W$ 
that does not depend on the choice of Coxeter generators $S$. 
\end{proof}

We call a subgroup of $W$ of the form $N({\cal F})$ a {\it basic characteristic subgroup}. 

\begin{corollary} 
If $W$ is a finitely generated Coxeter group 
and $N({\cal F})$ is a basic characteristic subgroup of $W$, 
then the isomorphism type of $W({\cal F}) = W/N({\cal F})$ 
is an isomorphism invariant of $W$. 
\end{corollary}

\section{The Spherical Rank Two Invariant}  

B. M\"uhlherr \cite{Muhlherr} has announced 
a solution of the isomorphism problem for finitely 
generated Coxeter groups $W$ such that 
$W$ has no basic subgroups of rank greater than 2 
with respect to some set of Coxeter generators. 
By the Basic Matching Theorem, if $W$ 
has no basic subgroups of rank greater than 2 
with respect to some set of Coxeter generators, 
then $W$ has no basic subgroups of rank greater than 2 
with respect to every set of Coxeter generators.  
Therefore it makes sense to say that $W$ has no basic 
subgroups of rank greater than 2 without regard to 
a set of Coxeter generators. 

In this section we describe a characteristic subgroup $N_2$ of 
a finitely generated Coxeter group $W$ such that $W_2 = W/N_2$ 
is a Coxeter group with no basic subgroups of rank greater than 2 
and such that the isomorphism type of $W_2$ is an isomorphism 
invariant of $W$.

Let ${\bf X}_n$ be one of the finite irreducible Coxeter systems  
${\bf A}_n$, ${\bf B}_n$, ${\bf C}_n$,  
${\bf E}_6$, ${\bf E}_7$, ${\bf E}_8$, ${\bf F}_4$, ${\bf G}_3$, ${\bf G}_4$ 
of rank $n\geq 3$. 
We now define a characteristic subgroup $N({\bf X}_n)$ of ${\bf X}_n$ 
for each ${\bf X}_n$. 
Let $N({\bf X}_n)$ be the commutator subgroup of the Coxeter group ${\bf X}_n$ 
if $n\geq 5$ or if ${\bf X}_n =  {\bf A}_4, {\bf G}_3$, or ${\bf G}_4$. 

Let $a_1, a_2, a_3$ be the Coxeter generators of ${\bf A}_3$ indexed so that  
$m(a_1,a_2) = m(a_2,a_3) = 3$.  
Let $N({\bf A}_3)$ be the normal closure 
in the group ${\bf A}_3$ of the element $a_1a_3$. 
Then $N({\bf A}_3)$ is a characteristic subgroup of ${\bf A}_3$ 
characterized by the property that $N({\bf A}_3)$ is the unique normal 
subgroup of ${\bf A}_3$ such that ${\bf A}_3/N({\bf A}_3)$ is isomorphic 
to ${\bf A}_2$ according to Table 3 of Maxwell \cite{Maxwell}. 

Let $b_1,b_2,b_3,b_4$ be the Coxeter generators of ${\bf B}_4$ 
indexed so that 
$$m(b_1,b_4) = m(b_2,b_4) = m(b_3,b_4) = 3.$$
Let $N({\bf B}_4)$ be the normal closure 
in the group ${\bf B}_4$ of the elements $b_1b_2$ and $b_2b_3$. 
Then $N({\bf B}_4)$ is a characteristic subgroup of ${\bf B}_4$ 
characterized by the property that $N({\bf B}_4)$ is the unique normal 
subgroup of ${\bf B}_4$ such that ${\bf B}_4/N({\bf B}_4)$ is isomorphic 
to ${\bf A}_2$ according to Table 3 of Maxwell \cite{Maxwell}. 

Let $c_1, c_2, c_3$ be the Coxeter generators of ${\bf C}_3$ such that 
$m(c_1,c_2) = 3$ and $m(c_2,c_3) = 4$.  
Let $N({\bf C}_3)$ be the normal closure 
in the group ${\bf C}_3$ of the element $(c_2c_3)^2$. 
Then $N({\bf C}_3)$ is a characteristic subgroup of ${\bf C}_3$ 
characterized by the property that $N({\bf C}_3)$ is the unique normal 
subgroup of ${\bf C}_3$ such that ${\bf C}_3/N({\bf C}_3)$ is isomorphic 
to ${\bf A}_2 \times {\bf A}_1$ according to Table 3 of Maxwell \cite{Maxwell}. 

Let $c_1,c_2,c_3,c_4$ be the Coxeter generators of ${\bf C}_4$ 
indexed so that 
$$m(c_1,c_2) = m(c_2,c_3) = 3\quad\hbox{and}\quad m(c_3,c_4) = 4.$$
Let $N({\bf C}_4)$ be the normal closure 
in the group ${\bf C}_4$ of the element $c_1c_3$. 
Then $N({\bf C}_4)$ is a characteristic subgroup of ${\bf C}_4$ 
characterized by the property that $N({\bf C}_4)$ is the unique normal 
subgroup of ${\bf C}_4$ such that ${\bf C}_4/N({\bf C}_4)$ is isomorphic 
to ${\bf A}_2 \times {\bf A}_1$ according to Table 3 of Maxwell \cite{Maxwell}. 

Let $f_1,f_2,f_3,f_4$ be the Coxeter generators of ${\bf F}_4$ 
indexed so that 
$$m(f_1,f_2) = m(f_3,f_4) = 3\quad\hbox{and}\quad m(f_2,f_3) = 4.$$
Let $N({\bf F}_4)$ be the normal closure 
in the group ${\bf F}_4$ of the element $(f_2f_3)^2$. 
Then $N({\bf F}_4)$ is a characteristic subgroup of ${\bf F}_4$ 
characterized by the property that $N({\bf F}_4)$ is the unique normal 
subgroup of ${\bf F}_4$ such that ${\bf F}_4/N({\bf F}_4)$ is isomorphic 
to ${\bf A}_2 \times {\bf A}_2$ according to Table 3 of Maxwell \cite{Maxwell}. 

Let $(W,S)$ be a Coxeter system of finite rank. 
Let $N(W)$ be the normal closure in $W$ of the subgroups 
$N(\langle{B}\rangle)$ defined above for every base $B$ of $(W,S)$ 
of rank greater than 2. Let $W^{(2)} = W/N(W)$. 
Let $\eta: W \to W^{(2)}$ be the quotient homomorphism, 
and let $S^{(2)} = \eta(S)$. 

\begin{theorem}  
The pair $(W^{(2)}, S^{(2)})$ is a Coxeter system 
that can be obtained from $(W,S)$ be a finite series 
of elementary edge quotient operations. 
The group $N(W)$ is a characteristic subgroup of $W$ 
that does not depend on the choice of Coxeter generators $S$.
\end{theorem}
\noindent\begin{proof}
Quotienting out the group $N(\langle B\rangle)$ for 
each base $B$ of $(W,S)$ of rank greater than 2 
can be realized by elementary edge quotient operations. 
Therefore $(W^{(2)}, S^{(2)})$ is a Coxeter system 
that can be obtained from $(W,S)$ be a finite series 
of elementary edge quotient operations. 

Let $S'$ be another set of Coxeter generators of $W$. 
By the Basic Matching Theorem and the characteristic 
properties of the groups $N({\bf X}_n)$, 
the group $N(W)$ defined in terms of the generators $S$ 
is the same as the group $N(W)$ defined in terms of the generators $S'$.   
Thus $N({\cal F})$ is a characteristic subgroup of $W$ 
that does not depend on the choice of generators $S$. 
\end{proof}

\begin{corollary} 
If $W$ is a finitely generated Coxeter group,  
then the isomorphism type of $W^{(2)} = W/N(W)$ 
is an isomorphism invariant of $W$. 
\end{corollary}

It may happen that $(W^{(2)}, S^{(2)})$ has a base of rank greater than 2.  
To get a quotient system with no bases of rank greater than 2, 
we may have to quotient out $N(W^{(2)})$, and then perhaps repeat 
the above quotienting process several times. 
This leads to a finite nested sequence 
$$\{1\}= N^{(1)}(W) \subset N^{(2)}(W) \subset \cdots \subset N^{(\ell)}(W)$$ 
of characteristic subgroups of $W$ such that if $W^{(i)} = W/N^{(i)}(W)$ 
and if $\eta_i:W \to W^{(i)}$ is the quotient homomorphism, 
then 
$$N^{(i+1)}(W) = \eta_i^{-1}(N(W^{(i)}))$$ 
for each $i=1,\ldots,\ell-1$, 
and $W^{(\ell)}$ has no basic subgroups of rank greater than 2, 
and $\ell$ is as small as possible. 
We have that $W^{(i+1)} = (W^{(i)})^{(2)}$ for each $i=1,\ldots,\ell-1$. 
Therefore the isomorphism type of $W^{(i)}$ for each $i=1,\ldots,\ell$ 
is an isomorphism invariant of $W$. 
It follows from the Basic Matching Theorem that $\ell$ does not depend 
on a choice of Coxeter generators for $W$, and so $\ell$ 
is an isomorphism invariant of $W$. 
We call $\ell$ the {\it spherical rank 2 class} of $W$. 
We have $\ell\geq 1$ with $\ell = 1$ if and only if $W$ 
has no basic subgroups of rank greater than 2. 
Figure 2 shows the P-diagrams of 
a sequence $W^{(1)},\ldots,W^{(\ell)}$ with $\ell = 4$ 
for the Coxeter group $W = W^{(1)}$.  

Define $N_2 = N^{(\ell)}(W)$. 
Then $N_2$ is a characteristic subgroup of $W$ 
such that $W_2 = W/N_2$ has no basic subgroups of rank greater than 2.  
The isomorphism type of $W_2$ is an isomorphism invariant of $W$ 
which we call the {\it spherical rank 2 isomorphism invariant} of $W$. 

Let $\eta:W \to W_2$ be the quotient homomorphism, 
and let $S_2 = \eta(S)$.  
Then $(W_2,S_2)$ is a Coxeter system that can be obtained from $(W,S)$ 
by a finite series of elementary edge quotient operations.

\medskip

$$\mbox{
\setlength{\unitlength}{.8cm}
\begin{picture}(16,5)(0,0)
\thicklines
\put(.5,1){\circle*{.15}}
\put(.5,1){\line(1,0){2}}
\put(2.5,1){\circle*{.15}}
\put(.5,1){\line(1,-1){1}}
\put(2.5,1){\line(-1,-1){1}}
\put(1.5,0){\circle*{.15}}
\put(.5,1){\line(0,1){3}}
\put(.5,2.5){\circle*{.15}}
\put(.5,4){\circle*{.15}}
\put(.5,2.5){\line(1,0){2}}
\put(2.5,2.5){\circle*{.15}}
\put(.5,4){\line(1,0){2}}
\put(2.5,4){\circle*{.15}}
\put(2.5,1){\line(0,1){3}}
\put(.5,.15){3}
\put(2.15,.15){5}
\put(1.4,1.2){2}
\put(1.4,2.7){2}
\put(1.4,4.2){2}
\put(0,1.5){3}
\put(0,3){3}
\put(2.7,1.5){3}
\put(2.7,3){3}
\put(5,2.5){\line(0,1){1.5}}
\put(5,2.5){\circle*{.15}}
\put(5,4){\circle*{.15}}
\put(4.5,3){3}
\put(5,2.5){\line(1,0){2}}
\put(7,2.5){\circle*{.15}}
\put(5.9,2.7){2}
\put(7,2.5){\line(0,1){1.5}}
\put(7,4){\circle*{.15}}
\put(7.2,3){3}
\put(5,4){\line(1,0){2}}
\put(5.9,4.2){2}
\put(5,2.5){\line(1,-1){1}}
\put(6,1.5){\circle*{.15}}
\put(5,1.65){3}
\put(7,2.5){\line(-1,-1){1}}
\put(6.65,1.65){3}
\put(9.5,4){\line(1,0){2}}
\put(9.5,4){\circle*{.15}}
\put(11.5,4){\circle*{.15}}
\put(10.4,4.2){2}
\put(9.5,4){\line(1,-1){1}}
\put(10.5,3){\circle*{.15}}
\put(9.5,3.15){3}
\put(11.5,4){\line(-1,-1){1}}
\put(11.15,3.15){3}
\put(10.5,1.5){\line(0,1){1.5}}
\put(10.5,1.5){\circle*{.15}}
\put(10.7,2.1){3}
\put(14,1){\line(0,1){3}}
\put(14,1){\circle*{.15}}
\put(14,2.5){\circle*{.15}}
\put(14,4){\circle*{.15}}
\put(14.2,1.5){3}
\put(14.2,3){3}
\end{picture}}$$

\medskip

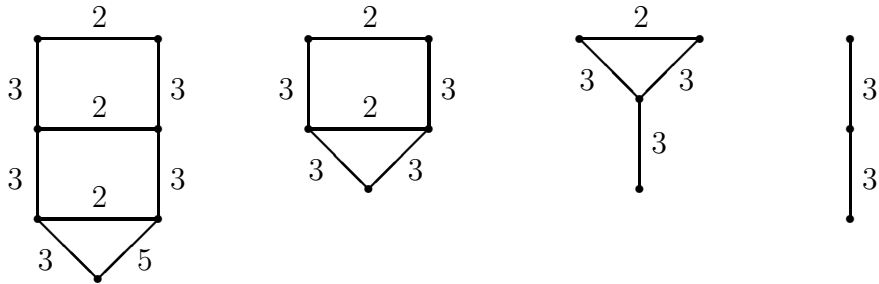
\centerline{\bf Figure 2}
\medskip

\newpage

\section{Conclusion}

Let $(W,S)$ be a Coxeter system of finite rank. 
In this paper, we have described three characteristic subgroups 
$N_b,N_{e},N_2$ of $W$ each leading to a quotient isomorphism invariant of $W$.  
It is interesting to note that 
$$N_2 \subset N_{e} \subset N_b,$$
and so the quotient isomorphism invariants corresponding to $N_b,N_{e},N_2$ 
are progressively stronger. 
The algorithm for finding a P-diagram for the system $(W_b,S_b)$ 
starting from a P-diagram of $(W,S)$ is computational fast. 

The algorithm for finding a P-diagram for the system $(W_{e},S_{e})$ 
is slower since it has to determine the bases of $(W,S)$ 
of type ${\bf D}_2(4q+2)$ that satisfy the conditions of Theorem 5.4;  
but, this algorithm is only slightly slower since the conditions 
of Theorem 5.4 are easy to check. 
The algorithm for finding a P-diagram for the system $(W_{e},S_{e})$ 
would most likely be incorporated into an efficient computer program 
that determines if two finite rank Coxeter systems have isomorphic groups, 
since the even isomorphism invariant would usually determine 
that two random finite rank Coxeter systems have nonisomorphic groups.

The algorithm for finding a P-diagram for the system $(W_2,S_2)$ 
is the slowest, but it is not much slower, since it only has to find a subdiagram  
of the P-diagram of $(W,S)$ of type ${\bf A}_3$, ${\bf C}_3$ or ${\bf G}_3$ 
before it performs an edge quotient operation 
on an edge of the subdiagram.  
If the subdiagram is of type ${\bf A}_3$ or ${\bf G}_3$, then 
the edge with label 2 is eliminated. 
If the subdiagram is of type ${\bf C}_3$, then the 4 edge label is reduced to 2. 
The algorithm then repeats the routine of searching for 
a subdiagram of type ${\bf A}_3$, ${\bf C}_3$ or ${\bf G}_3$ 
and performing the corresponding edge quotient operation. 

The algorithm for finding a P-diagram 
for the system $(W_2,S_2)$ would most likely be useful in an efficient program 
that determines if two finite rank Coxeter systems have isomorphic groups, 
since the solution of the isomorphism problem for finite rank 
Coxeter systems that have no bases of rank greater than 2 
is considerably simpler than any 
general solution of the isomorphism problem.

\end{document}